\documentclass[12pt]{article}
\makeatletter
\def\section{\@startsection{section}{1}{\z@}{-3.5ex plus -1ex minus -2.ex}
{2.3ex plus .2ex}{\Large\bf}}

\def\subsection{\@startsection{subsection}{2}{\z@}{-3.25ex plus
 -1ex minus -2.ex}
{1.5ex plus .2ex}{\bf}}

\def\vsn{\vskip 1pc \noindent}
\def\f{\newline}
\def\e{\varepsilon}

\def\rr{{\bf R}}
\def\nn{{\bf N}}

\newcommand{\be} {\begin{equation}}
\newcommand{\ee} {\end{equation}}
\newcommand{\bd} {\begin{displaymath}}
\newcommand{\ed} {\end{displaymath}}

\newcommand{\bq}{\begin{eqnarray}}
\newcommand{\eq}{\end{eqnarray}}
\newcommand{\bqn}{\begin{eqnarray*}}
\newcommand{\eqn}{\end{eqnarray*}}
\newcommand{\ba}[1]{\begin{array}{#1}}
\newcommand{\eqa}{\end{array}}
\def\qed{
   \\[-4ex]
  \hbox to \hsize{\hfill \vrule height 1.6ex width 1.5ex
  depth -.1ex}}

\begin{document}

\bibliographystyle{alpha}

\thispagestyle{empty}

\begin{center} {\Large {\bf
Adaptive mesh point selection for the efficient solution of scalar IVPs
 }\footnotemark[1] } 
\end{center}
\footnotetext[1]{ ~\noindent This research was partly supported by the Polish Ministry of Science 
and Higher Education \vsn}

\medskip
\begin{center}
{\large {\bf Boles\l aw Kacewicz \footnotemark[2]  }}
\end{center}
\footnotetext[2]{ 
\begin{minipage}[t]{16cm} 
 \noindent
{\it AGH University of Science 
and Technology, Faculty of Applied Mathematics,\\
\noindent  Al. Mickiewicza 30, paw. A3/A4, III p., 
pok. 301,\\
 30-059 Krakow, Poland 
\newline
E-mail: $\;\;$ kacewicz@agh.edu.pl }  
\end{minipage} }

\thispagestyle{empty} 
$~$
\vsn
\vsn
\begin{center} {\bf{\Large Abstract}} \end{center}
\noindent
We discuss  adaptive mesh point selection for the solution of scalar IVPs. We consider a method that is optimal 
in the sense of the speed of convergence, and aim at minimizing the local errors. Although the speed of convergence 
cannot be improved by using the adaptive mesh points 
compared to the equidistant points, we show that the factor in the error expression can be significantly reduced.
We obtain formulas specifying the gain achieved in terms of the number of discretization subintervals, as well as in terms of the 
prescribed level of the local error.  Both nonconstructive and constructive versions
of the adaptive mesh selection are shown, and a numerical example is given.
\newpage
\noindent
\section {\bf{\Large Introduction}} 
\noindent
We deal with the question  how much an  adaptive choice of mesh points pays off in the solution of initial-value problems
\be
z'(t)=f(z(t)),\;\; t\in [a,b],\;\;\;\; z(a)=\eta,
\label{1}
\ee
where $a<b$, $f:\rr \to \rr$  is a $C^r$ function and $\eta\in \rr$. 
\f
Numerical analysts have been using adaptive techniques in numerical codes for solving various problems with
considerable success.    Adaption is a standard tool in numerical packages, see e.g.  the package QUADPACK \cite{Pie} for numerical integration, or,
 among many others, the well known solver DIFSUB by C.W. Gear or the library ODEPACK by A. Hindmarsh devoted to various types of ODEs.
A measure of practical efficiency of an adaptive strategy  is usually  the  performance for a number of computational examples. 
A method is considered 'good' if it works well for  large number of problems, and fails in  small number of cases. Many papers 
have  reported  
the advantage of adaption over nonadaption in that sense,  to mention only as  samples the old paper  \cite{lyness},  
more recent one \cite{jackiewicz}, or \cite{mazzia}, where a mesh selection strategy is discussed for Runge-Kutta methods.
Such an approach obviously gives us  a considerable practical knowledge, but it is not complete.  
 Many step size control strategies are not supported by a theoretical analysis.
For instance, an important question how  a particular strategy
influences the cost of the process most often remains open. Recently, advantages of adaptive selection of mesh points
were rigorously studied for problems with singularities, see e.g.  \cite{KP}, \cite{PlaWas}.
\f
In this paper  we  present 
 results explaining  potential gain of   adaptive mesh point selection for a regular problem  (\ref{1}).
In particular, we  rigorously discuss the accuracy and cost of an adaptive process
for a well precised class of problems, not only for a number of computational examples.
For the  integration of scalar $C^4$ functions similar questions have recently been  addressed   for the Simpson rule in \cite{plaskota},
where it is shown  that the
adaptive mesh selection allows us to reduce the error by reducing the asymptotic constant of the method. 
Adaptive mesh points for the approximation of univariate $W^{2,\infty}$ functions is discussed in \cite{hickernell} .
\f
We consider in this work the $C^r$ right-hand side functions $f$ in (\ref{1}).  
It is well known that  in the worst-case or asymptotic settings, with $m+1$ mesh points   one can 
achieve  local errors  of order $m^{-(r+1)} $ as $m\to \infty$.   This  by standard means translates
to the global error $O(m^{-r})$. The exponent $-r$  has been  shown  best possible,
for details see e.g. \cite{JoC1}. 
Furthermore,  the best  speed of convergence as $m\to \infty$  can be achieved by
using the equidistant mesh points.  
We have to add that for systems of IVPs the {\it information} about $f$ that gives us the global error $O(m^{-r})$  must itself be adaptive, 
in spite of the fact that the mesh points can be equidistant, see \cite{JoC2} for  explanation of a difference between
adaptive mesh and adaptive information. 
\f
 In these results a constant in the $'O'$- notation depends 
on a class of functions $f$  in the worst-case setting, and on a particular $f$ in the asymptotic setting. 
The size of the constant is not controlled; it depends on  a {\it global}  behavior of derivatives of $f$ in the domain.
In order to reduce local errors, in the next sections we will include the constant in the $'O'$- notation to the analysis.
To study possible advantages of adaptive selection of mesh points,  we consider  one of the methods
with best convergence $O(m^{-r})$, given by (\ref{metoda}).  We show formulas for the local error of the
method, which will serve us to define mesh points with asymptotically  minimal maximum local error.  The selected points 
are adapted to a {\it local}  behavior of $f$.
We express the local errors  in terms of $m$, or, alternatively, for a given $\e>0$, we ask what  $m$ should be to 
achieve  local errors proportional to $\e$.
\f
The formulas obtained for the optimal mesh points  are not constructive. 
We next show  how the method can be modified to get computable mesh  points and approximations.  Compared to the 
'ideal' result, the local error bound and the cost bound are in this version increased by (known) factors dependent only on  $r$ (but not on $f$).
It turns out  that the adaptive choice of the mesh points allows us to achieve the maximum local error 
$\left((b-a)/m\right)^{r+1}\, S(m)$. The factor $S(m)$ is bounded from above and below by positive constants dependent on $f$,
so that it does not improve the rate of convergence. However,  the advantage of using adaptive mesh points is hidden in $S(m)$,
since the value of $S(m)$  can be much smaller for adaptive than nonadaptive points. 
\f
The paper is organized as follows. In Section 2 we define the class of functions $f$ and precise the aim of the paper. Section 3
presents the method under consideration and a convergence result. In Section 4 we give local error expressions which are
used in Section 5 to define optimal (nonconstructive) mesh points. Section 6 is devoted to a constructive modification
of the method which is finally described in the algorithm ADMESH. The error and cost properties of ADMESH are shown 
in Theorem 1 which summarizes the results of the paper. 
In Remark 3 we shortly comment on  generalization of the results to systems of IVPs.
The behavior of the algorithm ADMESH is illustrated in Section 7 by a numerical 
example. The experiment  shows how much one can gain using   adaptive mesh over the  equidistant mesh for a right hand 
side function with derivatives of varying magnitude in parts of the domain.  
\vsn
\section {\bf{\Large Problem formulation}} 
\noindent
Let $m\in \nn$.  We wish to compute approximations to the solution $z$ of (\ref{1}) 
at $m+1$ points $a=x_{0,m}< x_{1,m}<\ldots < x_{m,m}=b$, that is,  to find pairs $(x_{i,m},y_{i,m})$, $i=0,1,\ldots,m$, where
$y_{i,m}$ is a (good) approximation to $z(x_{i,m})$.  
Let $l(m)$ be any sequence convergent to $0$ as $m\to \infty$.
We consider for any $f$  a class of partitions of $[a,b]$. We assume that 
    there exist $K=K(f,a,b,\eta)$ and $k_0=k_0(f,a,b,\eta)$ such that for all $m\geq k_0$ and any partition it holds
\be
\max\limits_{0\leq i\leq m-1}(x_{i+1,m} - x_{i,m}) \leq K\, l(m).
\label{podzial}
\ee
Note that we always have $\max\limits_{0\leq i\leq m-1}(x_{i+1,m} - x_{i,m})\geq (b-a)/m$ for $m\geq 1$. Thus, the condition 
 (\ref{podzial}) implies that $l(m)$ cannot go to zero faster than $1/m$. The convergence of $l(m)$ can be arbitrarily slow, and 
the constant $K$ can be arbitrarily large.
\f
To  shorten the notation, we shall omit in the sequel the second
subscript $m$, remembering that the choice of  points $x_i$ and $y_i$ can be different for varying $m$.
\f
We denote by $z_i$  the solution of the local problem
\be
z_i'(t)=f(z_i(t)),\;\; t\in [x_i,x_{i+1}],\;\;\;\; z_i(x_i)=y_i.
\label{2}
\ee
If the pairs $(x_i,y_i)$ are outputs of a certain method, then 
the local errors of the method are given by $|z_i(x_{i+1})-y_{i+1}|$, $i=0,1,\ldots,m-1$. Our aim is to minimize 
the maximal local error 
\be
\max_{0\leq i\leq m-1}\, |z_i(x_{i+1})-y_{i+1}| \to {\rm min}
\label{3}
\ee
with respect to all possible choices of the mesh points $x_0,\ldots, x_m$, and to find minimizing (optimal) pairs $(x_i^*, y_i^*)$. 
\f
The class of right-hand functions $f$ under consideration is given as follows. For $r\in \nn$,
\be
F_r=\{f=1/g: \; g\in C^r(\rr), g \mbox{ and } g^{(r)} \mbox{ have constant sign in } \rr,  f \mbox{ is Lipschitz in } \rr
\}.
\label{class}
\ee
 We denote the Lipschitz constant of $f$ by $L$, and assume without loss of generality that
$f$ is a positive function. 
Regarding the constant sign of $g^{(r)}$, we note that the same assumption 
about constant sign of the fourth derivative ($r=4$) of the integrand was essential in \cite{plaskota}  in the
 analysis of adaptive integration of scalar $C^4$ functions .
\f 
In the next sections  we  aim  at choosing a subdivision of $[a,b]$, possibly adapting it to  a  local  behavior of $f$,
   in order to minimize the local errors.
Our goal   will be to propose a rigorous strategy of mesh point selection, keeping the cost 
 of the process under control, and to establish  possible gain of adaption.
\vsn
\section {\bf{\Large The method under consideration }} 
\noindent
We shall use the identity
\be
t-x_i = \int\limits_{y_i}^{z_i(t)} \frac{1}{f(y)}\, dy =  \int\limits_{y_i}^{z_i(t)} g(y)\, dy, \;\;\; t\in [x_i,x_{i+1}].
\label{identity}
\ee
For a positive $f$, the solutions $z$ and $z_i$ are increasing functions.  
\f
Let $r$ be even. 
For a given interval $[y_i,y_{i+1}]$, let 
$\hat g_i$ be  the Lagrange interpolation polynomial of degree $\leq r-2$ 
for $g$ in $[y_i,y_{i+1}]$ based on $r-1$ equidistant nodes $p_0=y_i, p_1,\ldots, p_{r-3},  p_{r-2}=y_{i+1}$  for $r\geq 4$,
and one node $p_0=(y_{i+1}+y_i)/2$ for $r=2$. For odd $r$, we define $p_0=y_i, p_1,\ldots, p_r=y_{i+1}$ as equidistant points in
$[y_i,y_{i+1}]$. The interpolation  polynomial $\hat g_i$ of degree $\leq r-1$ is now based on the nodes $p_0,\ldots, p_{r-1}$.
\f
Let $x_0=a$, $y_0=\eta$.  We shall study  
the following method  relating sequences $\{x_{i}\}$ and $\{y_{i}\}$~:  
\be
x_{i+1}-x_i= \int\limits_{y_i}^{y_{i+1}} \hat g_i(y)\, dy.
\label{metoda}
\ee
Note that the right-hand side is 
the Newton-Cotes type quadrature approximating ${\displaystyle \int\limits_{y_i}^{y_{i+1}}  g(y)\, dy}$, and it continuously depends
 on $y_{i+1}$.  
\f
We shall now derive convenient expressions  for the remainder of the Newton-Cotes formulas for even and odd $r$.
Denote $\hat e_i(y)= g(y)-\hat g_i(y)$.  For even $r$, $r\geq 4$, 
we recall that the remainder of the Newton-Cotes quadrature  is given by 
\be
\int\limits_{y_i}^{y_{i+1}}  \hat e_i(y)\, dy  = 
\frac{g^{(r)}(\xi_i)}{r!}  \int\limits_{y_i}^{y_{i+1}} (y-p_0)^2(y-p_1)\cdot\ldots \cdot (y-p_{r-2})\, dy ,
\label{NC}   
\ee
where $\xi_i \in [y_i,y_{i+1}]$.
Denoting $\Delta_i= y_{i+1}-y_i$ and changing variables $y=\Delta_i x+y_i$, $x\in [0,1]$, we get for $r\geq 4$
\be
\int\limits_{y_i}^{y_{i+1}}  \hat e_i(y)\, dy  =  \frac{g^{(r)}(\xi_i) }{ r!}  \Delta_i^{r+1} C_r,
\label{NC1}
\ee
where 
\be
C_r=  \int\limits_{0}^{1} (x-\bar p_0)^2\cdot\ldots \cdot (x-\bar p_{r-2})\, dx \;\;\;\;\; (C_r<0),
\label{Cr}
\ee
and $\bar p_j$ are equidistant nodes in $[0,1]$.  For $r=2$, (\ref{NC1}) holds with $C_2=1/12$.
\f
Let $r$ be odd, $r\geq 3$. We have for $y\in [y_i,y_{i+1}]$ 
$$
\hat e_i(y)= \frac{g^{(r)}(\xi_{i,y})}{r!} (y-p_0)\ldots (y-p_{r-1}),\;\;\;\;\mbox{ for some }\;\; \xi_{i,y}\in  [y_i,y_{i+1}].
$$
The second integral in the splitting 
$$
\int\limits_{y_i}^{y_{i+1}}  \hat e_i(y)\, dy  = \int\limits_{y_i}^{p_{r-1}} + \int\limits_{p_{r-1}}^{y_{i+1}} 
$$
 can be written as
$$
 \int\limits_{p_{r-1}}^{y_{i+1}} \hat e_i(y)\, dy = \frac{g^{(r)}(\eta_i)}{r!} 
\int\limits_{ p_{r-1}}^{y_{i+1}}    (y-p_0)\ldots (y-p_{r-1})\, dy ,\;\;\;\;\mbox{ for some }\;\; \eta_i\in  [y_i,y_{i+1}].
$$
Since 
$$
\int\limits_{y_i}^{p_{r-1}} (y-p_0)\ldots (y-p_{r-1}) \, dy =0,
$$
the first intergral is equal to
$$
 \int\limits_{y_i}^{p_{r-1}} \hat e_i(y)\, dy = 
\int\limits_{y_i}^{ p_{r-1}}  \frac{g^{(r)}(\xi_{i,y})}{r!}   (y-p_0)\ldots (y-p_{r-1})\, dy
$$
$$
= \frac{g^{(r)}(\eta_i)}{r!}     \int\limits_{y_i}^{ p_{r-1}}  \frac{g^{(r)}(\xi_{i,y}) - g^{(r)}(\eta_i) }{ g^{(r)}(\eta_i)}   (y-p_0)\ldots (y-p_{r-1})\, dy.
$$
We denote 
\be
\gamma_i^y= \frac{g^{(r)}(\xi_{i,y}) - g^{(r)}(\eta_i) }{ g^{(r)}(\eta_i)}  .
\label{pom44}
\ee
The quantity $\gamma_i^y$ continuously depends on $y$;  this follows  from the well known interpolation remainder formula
written with the use of the divided difference of $g$ expressed in the integral  form. 
Summing up the two integrals we get for odd $r$,  $r\geq 3$
\be
\int\limits_{y_i}^{y_{i+1}}  \hat e_i(y)\, dy = \frac{g^{(r)}(\eta_i)}{r!} 
\int\limits_{ p_{r-1}}^{y_{i+1}}    (y-p_0)\ldots (y-p_{r-1})\, dy \, (1+\kappa_i),
\label{441}
\ee
where
\be
\kappa_i= 
 \frac{     {\displaystyle \int\limits_{y_i}^{ p_{r-1}} \gamma_i^y    (y-p_0)\ldots (y-p_{r-1})\, dy } } {
{\displaystyle  \int\limits_{ p_{r-1}}^{y_{i+1}}    
(y-p_0)\ldots (y-p_{r-1})\, dy }  }.
\label{442}
\ee
It is easy to see that
\be
|\kappa_i|\leq \frac{2r^{r+1}}{(r-1)!} \sup\limits_{y\in [y_i,y_{i+1}]} |\gamma_i^y|.
\label{kap}
\ee
Changing variables as above in the integral in (\ref{441}), we get 
\be
\int\limits_{y_i}^{y_{i+1}}  \hat e_i(y)\, dy = \frac{g^{(r)}(\eta_i)}{r!} \Delta_i ^{r+1} C_r (1+\kappa_i),
\label{441a}
\ee
where 
\be
C_r= \int\limits_{ 1-1/r}^{1}    (x-\bar p_0)\ldots (x-\bar p_{r-1})\, dx \, ,
\label{441aa}
\ee
and $\bar p_j$, $j=0,1,\ldots, r$  are equidistant points in $[0,1]$.
It is easy to see that the case $r=1$ is included in (\ref{441a})  with  $\kappa_i=0$.
The formulas  (\ref{NC1})   and (\ref{441a}) will be used in the next section.
\vsn
For any  sequence $a=x_0<x_1<\ldots <x_{m}=b$ under consideration  there is a
  sequence $\eta=y_0<y_1<\ldots < y_{m}$  which satisfies (\ref{metoda}),
and has the global error bounded (as it can be expected) as in the following proposition.  
Some bounds obtained in the proof will be used in the sequel.
Let $h_i=x_{i+1}-x_i$, $i=0,1,\ldots,m-1$.
\vsn
{\bf Proposition 1}$\;\;$ {\it   Let $f\in F_r$. There is $m_0\in \nn$ such that for any $m\geq m_0$,
for any partition $a=x_0<x_1<\ldots <x_{m}=b$ satisfying (\ref{podzial})  there exists a sequence $\eta=y_0<y_1<\ldots < y_{m}$
satisfying (\ref{metoda}) such that   
\be
|y_i-z(x_i)|\leq M_i  \max\limits_{0\leq j\leq i-1} h_j^r,\;\;\;\; i=0,1,\ldots,m \;\;\;\; (\mbox{with } \max\limits_{0\leq j\leq  -1} = 1) .
\label{lemma1}
\ee
 Here $M_0=0$, and $M_{i+1}=  \exp(Lh_i) M_i +\tilde M h_i$, $i=0,1,\ldots, m-1$, where $\tilde M$ is given by (\ref{new1}).
Hence,   $M_i\leq M:=\exp(L(b-a))(b-a) \tilde M$. 
}
\vsn
{\bf Proof}$\;$ We prove (\ref{lemma1}) by induction with respect to  $i$. The statement holds for $i=0$. Suppose that there exist $y_0<y_1<\ldots < y_i$ satisfying (\ref{metoda})
and (\ref{lemma1}), and $M_i\leq M$.  Let 
$$
F(y) = \int\limits_{y_i}^y g(z)\, dz,\;\;\;\;  \hat F(y)=\int\limits_{y_i}^y \hat g_i(z,y)\, dz \;\;\;\; \mbox{  and }\;\;\;\; H(y)=F(y)-\hat F(y), \;\; y\geq y_i
$$
(the notation  $\hat g_i(\cdot,y)$  reflects the fact that the interpolation polynomial is defined on the interval $[y_i,y]$). 
Note that $F$ and $\hat F$ are continuous functions, $F'(y)=g(y)>0$, and
$F(z_i(x_{i+1}))=x_{i+1}-x_i$.  Our aim is to show the existence of a solution $y_{i+1}> y_i$ of the equation $\hat F(y) = x_{i+1}-x_i$.
Note that $\hat F(y_i) = 0< x_{i+1}-x_i$. We show that there is $\bar y>z_i(x_{i+1}) $ (which depends on $i$) such that
 $\hat F(\bar y) \geq  x_{i+1}-x_i$. This holds iff $H(\bar y)\leq F(\bar y) - F(z_i(x_{i+1}))$.  
Using (\ref{NC1}) or (\ref{441a}) 
with the interval $[y_i,y_{i+1}]$ replaced by $[y_i,y]$ we have that
\be
H(y)= \frac{ g^{(r)}(\bar\xi_{i,y}) } {r!} (y-y_i)^{r+1} C_r\, (1+\bar\kappa_i^y) \;\;\;\; (\mbox{where }  \bar\kappa_i^y=0 \mbox{ for even } r),
\label{Hy}
\ee
for some $ \bar\xi_{i,y}\in [y_i,y]$. Since $F( y) - F(z_i(x_{i+1}))= g(\tilde\xi_{i,y}) (y- z_i(x_{i+1}))$ for some $\tilde\xi_{i,y}\in {\rm conv} (z_i(x_{i+1}),y)$,
the equivalent condition on $\bar y$ reads
\be
\frac{ g^{(r)}(\bar\xi_{i,\bar y}) }{r!} (\bar y-y_i)^{r+1} C_r\, ( 1+\bar\kappa_i^{\bar y}) \leq  g(\tilde\xi_{i,\bar y}) (\bar y- z_i(x_{i+1})).
\label{bary}
\ee
We now show that (\ref{bary}) holds for $\bar y= y_i+2f(y_i)(x_{i+1}-x_i) $. The following auxilliary inequalities hold
for sufficiently large $m$  (where the starting value of $m$ only depends on$~$$f$)
\be
y_i\leq z(x_i)+|y_i- z(x_i)|\leq z(x_i)+M_i \max\limits_{0\leq j\leq i-1} h_j^r\leq z(b)+1 ,
\label{pom1}
\ee
$$
z_i(x_{i+1})\leq  z(x_{i+1})+|z_i(x_{i+1})- z(x_{i+1})| \leq z(x_{i+1})+\exp(Lh_i) |y_i-z(x_i)|
$$
\be
  \leq z(b)+1,
\label{pom2}
\ee
and
$$
\bar y = y_i + 2 f(y_i) h_i \leq z(x_i) + M_i \max\limits_{0\leq j\leq i-1} h_j^r + 2 f(y_i) h_i
$$
\be
\leq z(b)+1.
\label{pom4}
\ee
Let now 
$$
C=\sup\limits_{y\in [\eta, z(b)+1]} |g^{(r)} (y)| \;\;\; \mbox { and }\;\;\; c=\inf\limits_{y\in [\eta, z(b)+1]} g(y).
$$
We   come back to  (\ref{bary}).  In order to bound $|\bar\kappa_i^{\bar y}|$, we use (\ref{pom44}) and (\ref{kap}),
where the working variable $y$ in these formulas is replaced by $z$ and $y_{i+1}$ replaced by $y$,  with $y=\bar y$.
Since $g^{(r)}$ is uniformly continuous on $[\eta,z(b)+1]$, we have for $\bar y$  as above that 
                                      $|\bar\kappa_i^{\bar y}|\leq 1/2$ for $m$ sufficiently large,
where the limit value of $m$ only depends on $g$. 
A sufficient condition for (\ref{bary}) can now be written as
$$
C|C_r|\, 2^{r+1}/(c^{r+1}r!) (x_{i+1} -x_i)^r  \leq 1/3,
$$
  which holds true for sufficiently large $m$.
Consequently, there exists $y_{i+1}\in (y_i, y_i+2f(y_i)(x_{i+1}-x_i)  ]$ such that $\hat F(y_{i+1}) = x_{i+1}-x_i$, as claimed.
\f
It remains to show that
\be
|y_{i+1}-z(x_{i+1})| \leq M_{i+1}  \max\limits_{0\leq j\leq i} h_j^r
\label{indukcja}
\ee
and $M_{i+1}\leq M$.  We remember that 
$$
H(y_{i+1}) = F(y_{i+1}) - F(z_i(x_{i+1}))  =    g(\tilde\xi_{i,y_{i+1}})  (y_{i+1}-z_i(x_{i+1})),
$$
and $ \tilde\xi_{i,y_{i+1}} \in {\rm conv} (z_i(x_{i+1}), y_{i+1} )\subset [\eta,z(b)+1]$. From this 
the local error can be expressed as
\be
y_{i+1} - z_i(x_{i+1}) =  \frac{ g^{(r)}(\bar\xi_{i,y_{i+1}})} {r!}  \frac{1}{ g(\tilde\xi_{i,y_{i+1}}) }\,  (y_{i+1}-y_i)^{r+1} C_r \, (1+\bar\kappa_i^{y_{i+1}}).
\label{lokalny}
\ee
Taking into account that $y_{i+1}-y_i \leq 2f(y_i)(x_{i+1}-x_i)$  we get for sufficiently large $m$ that
\be
|y_{i+1} - z_i(x_{i+1})| \leq  \frac{C\,|C_r|} {cr!}\, (y_{i+1}-y_i)^{r+1} (3/2) \leq \tilde M h_i^{r+1},
\label{lokalny1}
\ee
where 
\be
\tilde M =  (3/2)C\,|C_r|\,2^{r+1}/(c^{r+2}r!). 
\label{new1}
\ee
Finally, the bound (\ref{indukcja}) on the global error
together with the formula for $M_{i+1}$ follow from the inductive assumption and the inequality
$$
|y_{i+1}-z(x_{i+1})| \leq |y_{i+1} - z_i(x_{i+1})| + |z_i(x_{i+1}) - z(x_{i+1})| \leq \tilde M h_i^{r+1} + \exp(Lh_i) |y_i - z(x_i)|.
$$
This holds for $m\geq m_0$, where $m_0$ only depends on $f$ (and is indeiendent of $i$).
The proof  is completed. \qed
\vsn
In particular, for $h_j=O(m^{-1})$ the global error of the method is
$$
\max\limits_{0\leq i\leq m} |y_i-z(x_i)| = O(m^{-r}), \;\;\; m\to \infty, 
$$
which is known to be optimal,  as far as the speed of convergence is concerned. The constant in the $'O'$-notation however depends on 
a global behavior of $f$, and it can be large, see the constant $\tilde M$ in the statement of Proposition 1. 
We now take into account 
a local behavior of $f$ in order to adjust the step sizes $h_j$ to the size of derivatives of $f$ in particular subintervals.
\vsn
\section {\bf{\Large Local error expressions }} 
\noindent
The local error of the method can be expressed due to (\ref{lokalny})  as 
\be
y_{i+1} - z_i(x_{i+1}) =   \frac{ g^{(r)}(\bar\xi_{i,y_{i+1}}) }{r!}  \frac{1}{ g(\tilde\xi_{i,y_{i+1}}) } \,  (y_{i+1}-y_i)^{r+1} C_r \, (1+\bar\kappa_i^{y_{i+1}})
\label{wyr1}
\ee
with $\bar\xi_{i,y_{i+1}} \in [y_i,y_{i+1}]$ and $\tilde\xi_{i,y_{i+1}} \in {\rm conv}(z_i(x_{i+1}), y_{i+1})$.  
\f
We shall adopt in what follows a convenient notation for relative errors used in the round off error analysis of numerical algorithms.
We have that $y_{i+1}-z_i(x_{i+1})= y_{i+1}-y_i + y_i - z_i(x_{i+1}) =   y_{i+1}-y_i      - f(z_i(\alpha_i)) h_i$, so that
$$
y_{i+1}-y_i= f(z_i(\alpha_i)) h_i ( 1+ \kappa_i),
$$
for some $\alpha_i\in [x_i,x_{i+1}]$, where  $\kappa_i= (y_{i+1}-z_i(x_{i+1}))/ (f(z_i(\alpha_i))\, h_i)$.
We  can alternatively write the local error as 
\be
y_{i+1} - z_i(x_{i+1}) =
\frac{ g^{(r)}(\bar\xi_{i,y_{i+1}}) }{r!   g(\tilde\xi_{i,y_{i+1}})\, ( g(z_i(\alpha_i)))^{r+1}   }\, C_r\,  
(x_{i+1}-x_i)^{r+1}\,(1+\bar\kappa_i^{y_{i+1}}) ( 1+ \kappa_i)^{r+1},
\label{wyr2}
\ee
where $z_i(\alpha_i)\in [y_i,z_i(x_{i+1})]$ and $\max\limits_{0\leq i\leq m-1}|\kappa_i|$ tends to zero as $m\to \infty$. 
The last convergence is uniform with respect to the class of partitions $\{x_i\}$.
\vsn
The following remarks  will be used in what follows.
\vsn
{\bf Remarks}$\;\;$ 
\f
{\bf 1.}  Let $\gamma :[\eta,z(b)+1] \to \rr$ be a continuous function of constant sign, and $[\alpha_{i,m},\beta_{i,m}]\subset [\eta,z(b)+1]$, 
$\alpha_{i,m}<\beta_{i,m}$, 
$i=0,1,\ldots, m-1.$  Assume that 
$\max\limits_{0\leq i\leq m-1} (\beta_{i,m}-\alpha_{i,m})$ tends to zero as $m\to \infty$.  Then by the uniform continuity,
for any $ z_1,z_2\in [\alpha_{i,m},\beta_{i,m}]$  we have that 
$$\gamma(z_1)=\gamma(z_2) (1+\bar\kappa_{i,m}),$$
 for some $\bar\kappa_{i,m}$, where
$\max\limits_{0\leq i\leq m-1} |\bar\kappa_{i,m}|$ tends to zero as $m\to \infty.$
\f
{\bf 2.}   Let 
$\lim\limits_{m\to \infty}\max\limits_{0\leq i\leq m-1} |\bar\kappa_{i,m}^j|=0$ for $j=1,2$. Define $\bar\kappa_{i,m}^3$ by
$$
1+\bar\kappa_{i,m}^3=(1+   \bar\kappa_{i,m}^1)(1+\bar\kappa_{i,m}^2) \; \mbox{ or }\;  1+   \bar\kappa_{i,m}^3 =(1+   \bar\kappa_{i,m}^1)/
(1+   \bar\kappa_{i,m}^2)\;  \mbox{ or }\;  1+   \bar\kappa_{i,m}^3 = (1+   \bar\kappa_{i,m}^1)^{r+1}.
$$
Then obviously $\lim\limits_{m\to \infty}\max\limits_{0\leq i\leq m-1} |\bar\kappa_{i,m}^3|=0$. \qed
\vsn
By these remarks we have the following lemma.
\vsn
{\bf Lemma 1}$\;\;$ {\it  The absolute  local error of the method (\ref{metoda}) is given by
\be
|y_{i+1} - z_i(x_{i+1})| =   c_i\, \Delta_i^{r+1} (1+\kappa_i^1) ,  
\label{lemma11}
\ee
where  $c_i= \sup\limits_{y\in [y_i,y_{i+1}]} \left( |g^{(r)}(y)|  / g(y)\right) \, |C_r|\,/r! $, and 
$\lim\limits_{m\to \infty}\max\limits_{0\leq i\leq m-1} |\kappa_i^1|=0$.
\f
Alternatively,
\be
|y_{i+1} - z_i(x_{i+1})| =   \bar c_i\,h_i^{r+1}  (1+\kappa_i^2) ,  
\label{lemma12}
\ee
where
$\bar c_i= \sup\limits_{y\in [z(x_i) ,z(x_{i+1}) ]}
\left( |g^{(r)}(y)|  / \left(  g(y)^{r+2}  \right) \right)\, |C_r|\,/r!$, and $\lim\limits_{m\to \infty}\max\limits_{0\leq i\leq m-1} |\kappa_i^2|=0$.
\f
The convergence of $\max\limits_{0\leq i\leq m-1} |\kappa_i^1|$ and $ \max\limits_{0\leq i\leq m-1} |\kappa_i^2|$ in (\ref{lemma11})
and (\ref{lemma12}), respectively, is uniform with respect to the partition.
 }
\vsn
{\bf Proof}$\;\;$ To show (\ref{lemma11}) we use (\ref{wyr1}). Let $s_i\in [y_i,y_{i+1}]$ be a point at which the supremum
in the definition of $c_i$ is achieved. We use Remark 1 with $\alpha_{i,m}=\min\{ y_i,z(x_i)\}$ and 
$\beta_{i,m}=\max \{ y_{i+1}, z_i(x_{i+1}),  z(x_{i+1})\}$, with  function    $\gamma(y)=|g^{(r)}(y)|$ or   $\gamma(y)=g(y)$
and a  point  $z_1$ suitably chosen, 
and  with $z_2$ fixed to be $z_2=s_i$.
The number  $\kappa_i^1$ absorbes all  numbers $\bar\kappa_{i,m}$ that appear when applying Remark 1, 
in  accordance with Remark 2. 
\f
To show (\ref{lemma12}), we use (\ref{wyr2}) and Remarks 1 and 2 in a similar way. \qed
\vsn
The unknown numbers $c_i$ and $\bar c_i$ depend on a local behavior of the function $g$.
\vsn
\section {\bf{\Large Adaptive (nonconstructive) selection of mesh points }} 
\noindent
We now show how to (approximately) minimize the maximal absolute local error  skipping for a moment the question whether 
the mesh points can be constructed or not. From (\ref{lemma12}) 
\be
|y_{i+1} - z_i(x_{i+1})| =   \bar c_i\,(x_{i+1}-x_i)^{r+1}  (1+\kappa_i^2) . 
\label{powt}
\ee
with 
$\bar c_i= \sup\limits_{y\in [z(x_i) ,z(x_{i+1}) ]} \left( |g^{(r)}(y)|  / \left(  g(y)^{r+2}  \right) \right)\, |C_r|/r!$. 
\f
We note that  for any $f$ and $\alpha\in (0,1/2)$ there exists $m_0$ such that for any $m\geq m_0$ 
and any partition $\{x_i\}$ under consideration it holds
\be
(1-\alpha) \max\limits_{0\leq i\leq m-1} \bar c_i (x_{i+1}-x_i)^{r+1} \leq \max\limits_{0\leq i\leq m-1} |y_{i+1}-z_i(x_{i+1})| 
\leq (1+\alpha) \max\limits_{0\leq i\leq m-1} \bar c_i (x_{i+1}-x_i)^{r+1} .
\label{prop21}
\ee
Consider now the minimization problem 
\be
\max\limits_{0\leq i\leq m-1} \bar c_i\, (x_{i+1} - x_i)^{r+1} \to \, {\rm MIN} \mbox{ with respect to } x_1,x_2,\ldots, x_{m-1} ,
\label{minimization1}
\ee
with $x_0=a$ and $x_m=b$.
\f
Define the functions of variables $x_1, x_2,\ldots, x_{m-1}$ by
$$
p(x_i,x_{i+1}) = \bar c_i  (x_{i+1} -x_i)^{r+1}, \;\;\; x_{i+1}\geq x_i,
$$
and
 $$
P_m(x_0, x_1,\ldots, x_m) =  \max\limits_{0\leq i\leq m-1} p(x_i,x_{i+1}) .
$$
Note that $p$ is a continuous function of $(x_i,x_{i+1})$, 
it is an increasing function of $x_{i+1}$ for fixed $x_i$, 
and a decreasing function of $x_i$ for fixed $x_{i+1}$.
The function $P_m$ is continuous on the compact set $a=x_0\leq x_1\leq\ldots \leq x_m=b$, so that
it attains its infimum for some $a=x_0^*\leq x_1^*\leq\ldots \leq x_m^*=b$. 
The corresponding $\bar c_i$ are denoted by   $\bar c_i^*$.   We note that the infimum  
$$
\inf\limits_{x_0,x_1,\ldots,x_m} P_m(x_0,x_1,\ldots,x_m)
$$
is a nonincreasing function of $m$, since $P_m(x_0,x_1,\ldots,x_m)$ is equal to $P_{m+1}(x_0,x_1,\ldots,x_m,x_m)$
for any $x_0\leq x_1\leq\ldots \leq x_m$. 
\vsn
{\bf Proposition 2}$\;\;$ {\it 
It holds 
\be
p(x_i^*,x_{i+1}^*)= \bar c_i^*\, (x_{i+1}^*   - x_i^*)^{r+1} = k_m^* = {\rm const}, \;\;\; i=0,1,\ldots, m-1.
\label{lemma31}
\ee
The number  $k_m^*$ equals the minimal value in (\ref{minimization1}), 
$
k_m^*= \min\limits_{x_0,x_1,\ldots,x_m} \max\limits_{0\leq i\leq m-1} \bar c_i\, (x_{i+1} - x_i)^{r+1}.
$
\f
The points $a=x_0^*, x_1^*,\ldots, x_{m-1}^*, x_m^*=b$ are unique.
\f
Furthermore,
\be
x_{i+1}^*-x_i^*= (b-a) \frac{ (1/\bar c_i^*)^{1/(r+1)} }{ \sum\limits_{i=0}^{m-1} (1/\bar c_i^*)^{1/(r+1)} } \;\;\;\; \mbox{ and }\;\;\;\;
k_m^*=\frac{ (b-a)^{r+1} } { \left( \sum\limits_{i=0}^{m-1} (1/\bar c_i^*)^{1/(r+1)} \right)^{r+1} }.
\label{lemma32}
\ee
}
\f
{\bf Proof}$\;\;$ The proof of (\ref{lemma31}) follows from the following observation.
Suppose that there is $i$ such that 
$$p(x_i^*,x_{i+1}^*) < p(x_{i+1}^*,x_{i+2}^*)$$
(the case $'>'$ is analogous).  Then 
we can decrease   $\max \{p(x_i^*,x_{i+1}^*), p(x_{i+1}^*,x_{i+2}^*)  \}$  by slightly  increasing $x_{i+1}^*$.
Applying this observation if necessary a number of times, we can also decrease  $P_m(x_0^*, x_1^*,\ldots, x_m^*)$,
which is a contradiction.
\f
Given $x_i^*$, the point $x_{i+1}^*$ is a solution of $p(x_{i}^*,x_{i+1})=k_m^*$, where 
$k_m^*=\inf\limits_{x_0,x_1,\ldots,x_m} P_m(x_0, x_1,\ldots, x_m) $. The solution is unique, since 
$p(x_i^*, \cdot)$ is an increasing function.
\f
Relations (\ref{lemma32}) follow from (\ref{lemma31}) and the fact that $\sum\limits_{i=0}^{m-1}(x_{i+1}^* - x_i^*) = b-a.$ 
\f
\qed
\vsn
The sequence $\{k_m^*\}$ is nonincreasing.  A convenient expression for $k_m^*$ is the following
\be
k_m^*=\left( \frac{ b-a }{m} \right)^{r+1}  S(m) , 
\label{k-m}
\ee
where
\be
S(m)=
{\displaystyle
\frac {1}{ \left(  \frac{1}{m}  \sum\limits_{i=0}^{m-1} (1/\bar c_i^*)^{1/(r+1)} \right)^{r+1} }.
\label{k-m1}
}
\ee
Note that $S(m)$ plays here the role of a constant, since the dependence on $m$ is weak: for any $m\geq 1$
we have 
\be
0< c(f)\leq S(m)\leq C(f),
\label{Sm}
\ee
where 
$$
c(f)= \inf\limits_{y\in [\eta,z(b)]} b_g(y) \;\; \mbox{ and }\;\;    C(f)=  \sup\limits_{y\in [\eta,z(b)]} b_g(y)
$$
with $b_g(y)= \left( |g^{(r)}(y)|  / \left(  g(y)^{r+2}  \right) \right) |C_r|/r!$.
We see that the factor $S(m)$ in (\ref{k-m}) 
 does not  improve the speed of convergence
of $k_m^*$, which 
remains of order $\Theta\left(m^{-(r+1)} \right)$ as $m\to \infty$. However, the gain in the coefficient can be significant.
The number  $C(f)$ in the 'a priori'  bound  (\ref{Sm}) 
 can be large;
 it reflects a global behavior of the function $g$
in the entire interval $[\eta,z(b)]$.
This bound is sharp if $\bar c_i^*$ are essentially constant. In the opposite case, when $\bar c_i^*$ are all small except 
for a single large 
one equal to $C(f)$, the sum $\sum\limits_{i=0}^{m-1} (1/\bar c_i^*)^{1/(r+1)}$ can be much larger than
$$ 
 m \left(\frac{ 1} {\max\limits_{0\leq i\leq m-1} \bar c_i^*} \right)^{1/(r+1)} .
$$
 In this case $S(m)$ is much smaller than $C(f)$ and   $k_m^*$ is much smaller than the 'a priori'  bound,  
\be
k_m^* <<   \left( \frac{b-a}{m} \right)^{r+1} C(f) .
\label{111}
\ee
The gain from adjusting the mesh points to $\bar c_i^*$ is then   significant. 
\f
For comparison, consider the equidistant mesh points, $x_i=a+i(b-a)/m$. Then     
$$
\max\limits_{0\leq i\leq m-1} \bar c_i\, (x_{i+1} - x_i)^{r+1}= 
\max\limits_{0\leq i\leq m-1} \bar c_i \left( \frac{b-a}{m}\right)^{r+1}= C(f) \left( \frac{b-a}{m}\right)^{r+1}.
$$
Hence, for the equidistant mesh the 'a priori' upper bound in (\ref{111}) is attained. The points $x_i^*$ defined
in Proposition 2 can do much better.
\f
We stress again that the~ points $x_i^*$  for which the gain is achieved depend on 
unknown quantities; we do not show at this point how to construct them.  
\f
Note also that 
$$
\max\limits_{0\leq i\leq m-1} (x_{i+1}^*-x_i^*)   \leq  \frac{b-a}{m}  \cdot \left(\frac{ C(f)}{c(f)}\right)^{1/(r+1)} ,
$$
so that $\{x_i^*\}$ is an admissible partition (for any  $l(m)\geq 1/m$ and $K\geq (b-a)\left(C(f)/c(f)\right)^{1/(r+1)}$ in (\ref{podzial})).
\f
In many cases we are interested in computing approximations with the absolute local error not 
exceeding  a prescribed level $\e\in (0,1)$. 
That is, we wish 
to find the minimal number $m=m(\e)$ such that $k_m^*\leq \e$.  Hence, 
$$
k_m^*=\left( \frac{ b-a }{m} \right)^{r+1}  S(m)\leq \e,
$$
which gives us  that $m(\e)$ is the minimal $m$ such that
\be
m \geq   (b-a)S(m)^{1/(r+1)} \left( \frac{1}{\e} \right)^{1/(r+1)} .
\label{meps}
\ee
The 'a priori' bounds on $S(m)$ lead to 'a priori' bounds on $m(\e)$
\be
(b-a) c(f)^{1/(r+1)}  \left( \frac{1}{\e} \right)^{1/(r+1)} \leq m(\e) < (b-a) C(f)^{1/(r+1)} \left( \frac{1}{\e} \right)^{1/(r+1)} +1,
\label{meps1}
\ee
so that $m(\e)=\Theta\left( \left( 1/\e\right)^{1/(r+1)}\right)$ as $\e\to 0$. The actual value of $m(\e)$ can be 
however much
smaller than the upper bound, since $S(m)$ for all $m$ can be much smaller than $C(f)$.
\f
It is clear that the number of subdivision intervals $m(\e)$ will be crucial  for establishing the minimal cost of computing a
constructive  approximation with the absolute  local error at most $\e$.
\f
Proposition 2 leads to the following result about  minimization of the maximal absolute local error. Let
\be
L^m= \min\limits_{x_0,x_1,\ldots,x_m} \max\limits_{0\leq i\leq  m-1} |y_{i+1}-z_i(x_{i+1})|.
\label{minlokal}
\ee
The value of $L^m$ is asymptotically equal  to $k_m^*$, up to an arbitrarily small positive constant $\alpha$.
\vsn
{\bf Proposition 3}$\;\;$ {\it For any $f$ and $\alpha\in (0,1/2)$ there exists $ m_0$ such that for any $m\geq  m_0$ 
 the minimal error satisfies 
\be
 (1-\alpha) k_m^* \leq  L^m    \leq (1+\alpha) k_m^* ,
\label{prop22}
\ee
and 
\be
(1-2\alpha) \max\limits_{0\leq i\leq  m-1} |y_{i+1}^*-z_i(x_{i+1}^*)|  \leq  L^m   \leq 
\max\limits_{0\leq i\leq m-1} |y_{i+1}^*-z_i(x_{i+1}^*)|.
\label{prop23}
\ee
Hence, up to a (possibly small) constant $\alpha$, the mesh points $\{x_i^*\}$ are optimal.
\f
(Here $\{y_i^*\}$ are given  for $\{x_i^*\}$ by (\ref{metoda}), and 
 $z_i$ denotes the solution of the local problem with the initial condition $z_i(x_i^*)=y_i^*$.)
}
\vsn
{\bf Proof}$\;\;$  The proof follows from   (\ref{prop21}).  \qed
\vsn
Hence, the quantity $k_m^*$ is equal to the minimal maximum local error,  and the points $x_i^*$ define
the best partition  (up to the constant $\alpha$). 
The method (\ref{metoda}) needs at each step to compute
 the interpolation polynomial $\hat g_i$. 
In order to have local errors at level $\e$, 
the cost is thus at least $r\, m(\e)$ evaluations of the function $g$.
In the next section we  effectively construct the mesh points and modify (\ref{metoda}) to  compute approximations with the
local errors proportional to $\e$, with cost proportional to $m(\e)$.
\vsn
\section {\bf{\Large Adaptive constructive selection of mesh points }} 
\noindent
Let $\e\in (0,1)$. We shall slightly modify (\ref{metoda}) by replacing the interpolation polynomial 
$\hat g_i$ (defined on $[y_i,y_{i+1}]$) by an interpolation
 polynomial $\hat g_i^1$ defined on the interval dependent only on $y_i$. The approximation to $y_{i+1}$ will be obtained by
a number of steps of the bisection method. The replacement will allow us to use the same polynomial in all iterations.
This makes it possible to  avoid the $\log 1/\e$ factor in the cost bound, at expence of an additional factor dependent only on
$r$ in the error bound. 
The dependence of the cost on $g$  and $\e$, and  possible gain discussed in the previous section,  will be hidden in 
the quantity $m(\e)$. 
\f
To be specific,  we define points $\hat x_i$  as follows.  
 We set $\hat x_0=a$, $\hat y_0=\eta$. For a given $\hat x_i$ and $\hat y_i$, we compute the divided difference 
$$
g[\bar z_0^i,\bar z_1^i,\ldots,\bar z_r^i] ,
$$
where $\bar z_j^i$ are equidistant points from $[\hat y_i, \hat y_i + \e^{1/(r+1)}]$ (including the end points).
\f
Then we set
\be
\hat c_i=\hat c_i(\hat y_i) = \frac{2^{r+1} |g[\bar z_0^i,\bar z_1^i,\ldots,\bar z_r^i]| } {g(\hat y_i)^{r+2}}.
\label{cc}
\ee
The point $\hat x_{i+1}$ is defined as  the solution of
\be
\hat c_i (\hat y_i) (\hat x_{i+1} -\hat x_i)^{r+1} =   \frac{ 2^{r+1}}{ |C_r| }\, \frac{1}{1-\alpha}\, \e,          \;\;\; i=0,1,\ldots .
\label{224}
\ee
 Let  $\bar y_i= \hat y_i +2f(\hat y_i) (\hat x_{i+1} -\hat x_{i}) $.
Let $\hat g^1_i$ be the Lagrange interpolation polynomial for $g$ of degree at most $r-1$ based on $r$ 
equidistant points from  $[\hat y_i, \bar y_i]$ for $r\geq 2$ (including the end points), and  $\hat g^1_i(y)\equiv g(\hat y_i)$ for $r=1$.
We define  $\hat y_{i+1}$ as the solution of 
$$
\hat F^1 (y):= \int\limits_{\hat y_i}^y \hat g^1_i(z)\, dz = \hat x_{i+1} -\hat x_i 
$$
in the interval $[\hat y_i, \bar y_i ]$.  The  existence of $\hat y_{i+1}$ follows (for sufficiently small $\e$) 
from the arguments used in the proof of
Proposition 1 with $\{x_i\}$ and $\{y_i\}$ replaced by $\{\hat x_i\}$ and $\{\hat y_i\}$. We use the function 
 $H^1(y)= F(y)-\hat F^1(y)$, where 
$$ F(y)=\int\limits_{\hat y_i}^y g(z)\, dz.$$
By the standard interpolation error formula we have 
$$
|H^1(y)|\leq     \sup\limits_{\xi\in    [\hat y_i,\bar y_i]} |g^{(r)} (\xi)| (\bar y_i -\hat y_i)^r (y- \hat y_i)  /r! 
$$
$$
    =  \sup\limits_{\xi\in    [\hat y_i,\bar y_i]} |g^{(r)} (\xi)| 2^r (\hat x_{i+1}-\hat x_i)^r f(\hat y_i)^r (y- \hat y_i)  /r!   , \;\;\;     y\in [\hat y_i,\bar y_i].
$$
We define $\hat m$ to be the minimal $i$ for which $\hat x_i\geq b.$
\f
Consider now the  local errors $\hat y_{i+1}-z_i(\hat x_{i+1})$ of the pairs $(\hat x_i, \hat y_i)$  
($z_i$ is the solution of (\ref{2}) such that $z_i(\hat x_i)=\hat y_i$). 
Similarly as in the proof of Proposition 1, we have that
$$
\hat y_{i+1}- z_i(\hat x_{i+1}) = H^1 (\hat y_{i+1})/ g(\eta_i),
$$
for some $\eta_i\in {\rm conv} (\hat y_{i+1}, z_i(\hat x_{i+1}) )$. Hence
$$
|\hat y_{i+1}- z_i(\hat x_{i+1}) | \leq \gamma_i (\hat x_{i+1}-\hat x_i)^{r+1},
$$
where
$$
\gamma_i =\frac{ \sup\limits_{\xi\in    [\hat y_i,\bar y_i]} |g^{(r)} (\xi)| \, f(\hat y_i)^{r+1} \, 2^{r+1} } 
{ r! \inf\limits_{\xi\in    {\rm conv}(\hat y_{i+1}, z_i(\hat x_{i+1})) } g (\xi) }.
$$
We now use Remarks 1 and 2 to replace $\gamma_i$ by $\hat c_i$ to get 
\be
|\hat y_{i+1}- z_i(\hat x_{i+1}) | \leq \hat c_i (\hat x_{i+1}-\hat x_i)^{r+1}  (1+\kappa_i^2) ,
\label{2225}
\ee
with some $\kappa_i^2$ such that $\max\limits_{0\leq i\leq \hat m -1} |\kappa_i^2|$ tends to $0$ as $\e\to 0$. 
This together with the definition of 
$\hat x_{i+1}$ yields that for sufficiently small $\e$  we have
\be
|\hat y_{i+1}- z_i(\hat x_{i+1}) | \leq \frac{1+\alpha}{1-\alpha}\, \frac{2^{r+1}}{|C_r|} \, \e . 
\label{2226}
\ee
We now show that $\hat m \leq m(\e)$, where 
 $m(\e)$, which is our reference quantity, has been defined to be the minimal value of $m$ such that $k_m^*\leq \e$.
This will follow from the fact that for the optimal points we have $x_i^*\leq  \hat x_i$, $i=0,1,\ldots .$ 
Indeed, this holds for $i=0$. Let  $x_i^*\leq  \hat x_i$ for some $i$. If $x_{i+1}^*\leq \hat x_i$, then obviously $x_{i+1}^*\leq \hat x_{i+1}$,
so that it suffices to consider the case  $ \hat x_i<x_{i+1}^*$. Then
$$
\e \geq k_m^*= \bar c_i(x_i^*,x_{i+1}^*) (x_{i+1}^* -x_i^*)^{r+1} \geq \bar c_i(\hat x_i,x_{i+1}^*) (x_{i+1}^* -\hat x_i)^{r+1}.
$$
 For convenience, we have explicitly written here the arguments that 
$\bar c_i$  depends on. 
It follows from Remarks 1 and 2 that  for sufficiently small $\e$
$$
\e \geq \bar c_i(\hat x_i, x_{i+1}^*) (x_{i+1}^*- \hat x_i)^{r+1} =   \hat c_i(\hat y_i) (x_{i+1}^* -\hat x_i)^{r+1} |C_r| (1+\kappa_i^3)/2^{r+1}
$$
$$
\geq \hat c_i(\hat y_i) (x_{i+1}^* -\hat x_i)^{r+1} |C_r| (1-\alpha)/2^{r+1}.
$$
Thus,
$$
\hat c_i(\hat y_i) (x_{i+1}^* -\hat x_i)^{r+1} \leq 2^{r+1}/( |C_r| (1-\alpha))\, \e,
$$
which yields that $\hat x_{i+1}\geq x_{i+1}^*$, as claimed.
\vsn
It remains to show how we compute an approximation to $\hat y_{i+1}$. We apply the bisection method to the
equation $\hat F^1(y)= \hat x_{i+1} - \hat x_i$, starting from the interval $[\hat y_i, \hat y_i+2f(\hat y_i) (\hat x_{i+1}-\hat x_i) ]$
(see the proof of Proposition 1). 
After $l_i\geq 1$ steps the length of the interval is reduced to $f(\hat y_i)(\hat x_{i+1}-\hat x_i)/2^{l_i-1}$. We choose $l_i$ to the minimal number
such that
\be
f(\hat y_i)(\hat x_{i+1}-\hat x_i)/2^{l_i-1} \leq \e/2.
\label{kroki}
\ee
Equivalently, inserting $\hat x_{i+1}-\hat x_i$, we get that $l_i$ is the minimal number such that
\be
\frac{ 8f(\hat y_i)} {\hat c_i^{1/(r+1)}\, |C_r|^{1/(r+1)} }\, \left(\frac{1}{1-\alpha} \right)^{1/(r+1)} \, \e^{1/(r+1)-1} \leq 2^{l_i}.
\label{kroki1}
\ee
We note that  the bisection process does not require any new evaluations of $g$.
\f
Any point of the last bisection interval can be taken as an approximation to $\hat y_{i+1}$. 
We denote the selected approximation by the same symbol $\hat y_{i+1}$, and get
\be
|\hat y_{i+1} - z_i(\hat x_{i+1})| \leq 
       \left( \frac{1+\alpha}{1-\alpha}\, \frac{2^{r+1}}{|C_r|} +\frac{1}{2}\right)\, \e .
\label{225a}
\ee
We will refer to the above bisection procedure  as BISEC.
These considerations are summarized in the following algorithm   for solving (\ref{1}).
\vsn
\newpage\noindent
Algorithm {\bf ADMESH}
\vsn
\begin{tabular}{ll}
1$\;\;\;\;$ & Set $\e\in (0,1)$, $\alpha\in (0,1/2)$, $\hat x_0=a$, $\hat y_0=\eta$, $i=-1$\\
2$\;\;\;\;$ & $i:=i+1$\\
3$\;\;\;\;$ & Compute $d_r^i= g[\bar z_0^i, \bar z_1^i, \ldots, \bar z_r^i]$, where  $\bar z_j^i$ are equidistant points 
\\ $~$ &  in $[\hat y_i, \hat y_i + \e^{1/(r+1)}  ]$ 
(including the end points)\\ 
4$\;\;\;\;$ & Compute 
$
\hat c_i= 2^{r+1} |d_r^i| f(\hat y_i)^{r+2}  
$
and 
$ \hat x_{i+1} =\hat x_i +    2/\left( |C_r| \, \hat c_i  \, (1-\alpha)\right)^{1/(r+1)} \, \e^{1/(r+1)}$.
\\
$~$&If $\hat x_{i+1}\geq b$ then $\hat x_{i+1}:= b$\\
5$\;\;\;\;$& Compute the interpolation polynomial $\hat g^1_i$
\\
6$\;\;\;\;$ &Compute $\hat y_{i+1}$ by the algorithm BISEC  applied to the equation $\hat F^1(y)= \hat x_{i+1} - \hat x_i$\\
$~$&with $l_i$ steps starting from $[\hat y_i, \hat y_i +2 f(\hat y_i)(\hat x_{i+1} - \hat x_i)]$, with $l_i$  given by (\ref{kroki1}). \\
$~$&If $\hat x_{i+1}= b$ then go to STOP\\
7$\;\;\;\;$ & Go to 2\\
STOP &$~$ 
\end{tabular}
\vsn
The following theorem summarizes the error and cost properties of the algorithm ADMESH.
\vsn
{\bf Theorem 1}$\;\;$ {\it Let   $f\in F_r$  and  $\alpha \in (0,1/2)$.
There exists $\e_0=\e_0(f, \alpha)$ such that for any $\e\leq \e_0$ 
the algorithm ADMESH  computes  pairs $(\hat x_i,\hat y_i)$, $i=0,1,\ldots, \hat m$,  with the following error/cost
properties.  The maximum local error is bounded by
\be
\label{thm1}
\max\limits_{0\leq i\leq \hat m -1} |\hat y_{i+1} - z_i(\hat x_{i+1})| 
\leq \left( \frac{1+\alpha}{1-\alpha}\, \frac{2^{r+1}}{|C_r|} +\frac{1}{2}\right)\, \e ,
\ee
 where $z_i$ is the solution of the local problem (\ref{2}) with the initial condition $z_i(\hat x_i)=\hat y_i$.
\f
The cost of the algorithm ${\rm cost}(f,\alpha,\e)$ measured by the number of evaluations of  $f$ is bounded by 
\be
{\rm cost}(f,\alpha,\e) \leq 2r\, m(\e),
\label{thm12}
\ee
where $m(\e)$ is the (almost optimal) number of subintervals given in (\ref{meps}).
}
\vsn
{\bf Proof}$\;\;$ The bound (\ref{thm1}) follows from (\ref{225a}). The cost related to the $i$th interval consists of $r+1$ function evaluations to
compute $d_r^i$ and additional $r-1$ function evaluations to compute the interpolation polynomial $\hat g^1_i$. This and the bound $\hat m\leq m(\e)$  yield the cost bound (\ref{thm12}).
\f \qed 
\vsn
We comment on this result. Note first that  $\alpha$ can be an arbitrary small positive number which does not play
any crucial role. The accuracy achieved by the algorithm ADMESH  differs from  the 
accuracy achieved by the almost optimal points $x_i^*$  only by the explicitly known  factor dependent on $r$ (and
independent of $f$), see Proposition 3. The cost of ADMESH is proportional, with coefficient $2r$, to the reference value $m(\e)$.
 It follows from 
the discussion after Proposition 2, see (\ref{meps1}),  that  the 'a priori' upper bound on the cost is
\be
{\rm cost}(f,\alpha,\e) \leq 
2r\left( (b-a) C(f)^{1/(r+1)} \left( \frac{1}{\e} \right)^{1/(r+1)} +1\right).
\label{wniosekkoszt}
\ee
This upper bound is essentially achieved by the equidistant mesh.  
The advantage of the mesh points constructed in the algorithm ADMESH lies in the fact  that $m(\e)$, where the dependence on $f$ is hidden, 
can be much
smaller than the upper bound given in (\ref{meps1}), see the discussion after Proposition 2. Consequently, the actual cost
of getting the accuracy proportional to $\e$ can be much smaller than the upper bound in  (\ref{wniosekkoszt}).
\f
Note also that we can have the error bound in (\ref{thm1}) equal to a given number $\e_1$, by running the algorithm
with $\displaystyle { \e:=\e_1\left( \frac{1+\alpha}{1-\alpha}\, \frac{2^{r+1}}{|C_r|} +\frac{1}{2}\right)^{-1} } $.
\vsn
{\bf Remark 3}$\;\;$ It would be   of interest to generalize the above results to systems of IVPs. One can see that 
a straightforward generalization is not possible, since  there is no counterpart of  (\ref{identity}) for systems of IVPs.
Preliminary analysis however indicates that a progress in that direction is possible using a different technique.
This will be a topic of our future work. 
\vsn
\section {\bf{\Large Numerical example}} 
\noindent
To illustrate the behavior of ADMESH, we consider a problem with $r=2$, dependent on a parameter $\delta>0$
\be
z'(t)=\frac{3}{4} (z(t)-1)^{-3/2}, \;\; t\in [0,1],\;\;\; z(0)=1+\delta.
\label{test}
\ee
The right hand side function $f$
has the form $f=1/g$, where $g(z)=(4/3)(z-1)^{3/2}$. The second derivative $g''(z)=(z-1)^{-1/2}$ taken at the initial condition 
grows to infinity  as $1/\sqrt{\delta}$ with $\delta \to 0^+$. For such a function and small $\delta$  we should observe a significant advantage of
adaptive mesh points over the equidistant points. The solution satisfying the initial condition $z(x)=y\;\;$ ($x\geq 0$, $y>1$)
is given by 
$$
z(t)= \left(\frac{15}{8} (t-x)+(y-1)^{5/2}\right)^{2/5}+1.
$$
The testing program was translated to the C$^{++}$ code by P. Morkisz and B. Bo\.zek.
\f
The following table shows results computed by ADMESH for number of values of $\e$ and $\delta$. In the successive columns we
show the values of IADAPT (the number of adaptive mesh points), MAXERR (the maximum local error), MAXERR/BOUND (BOUND is the upper bound given in Theorem 1),  MAXERRG (the maximum global error), EQUIDIST/MAXERR (EQUIDIST is the maximal local error 
obtained with  $2*$IADAPT equidistant mesh points),  EQUIDISTG/MAXERRG (the same ratio for the maximal global errors). Since ADMESH requires $4$ evaluations 
of $f$ in each subinterval, and the equidistant mesh algorithm $2$ evaluations, the results in the latter case are computed for twice as much points. The computer precision is 
$10^{-16}$. We took $\alpha=0.25$.
\vsn
\begin{tabular}{llllllll}  
$\e$ & $\delta$ &{\tiny IADAPT} &{\tiny MAXERR} & $\frac{\mbox{{\tiny MAXERR}}}{\mbox{{\tiny BOUND}}}$ &{\tiny MAXERRG} 
&$\frac{\mbox{{\tiny EQUIDIST}}}{\mbox{{\tiny MAXERR}}}$  
&$\frac{\mbox{{\tiny EQUIDISTG}}}{\mbox{{\tiny MAXERRG}}}$ \\$~$\\
 0.01& 0.1& 5& 0.23& 0.014& 0.035& 7.39& 4.9\\
0.01& $10^{-4}$&11& 0.02& 0.011& 0.032& 19.57& 11.96\\
0.01& $10^{-8}$& 11& 0.02& 0.012& 0.039& 17.79& 11.00\\
$10^{-4}$&0.1&15&7.3$*10^{-4}$&0.046& 3.1$*10^{-3}$& 90.56&21.06\\
$10^{-4}$&$10^{-4}$&27& 6.7$*10^{-4}$&0.042&3.0$*10^{-3}$& 369.89& 84\\
$10^{-4}$&$10^{-8}$&30&  6.7$*10^{-4}$&0.042&3.0$*10^{-3}$& 371.69& 101\\
$10^{-8}$&0.1&252& 1.09$*10^{-7}$& 0.068& 8.24$*10^{-6}$&8291& 109\\
$10^{-8}$&$10^{-4}$&418& 1.85$*10^{-7}$&0.115& 8.28$*10^{-6}$&436463&9732\\
$10^{-8}$&$10^{-8}$&435& 2.30$*10^{-7}$&0.143& 8.32$*10^{-6}$&373152&12562\\
$10^{-16}$& 0.1&115332& 1.59$*10^{-15}$& 0.099 & 3.87$*10^{-11}$& 17051 & 95\\
$10^{-16}$& $10^{-4}$&192546& 1.59$*10^{-15}$& 0.099& 3.90$*10^{-11}$& 3.7$*10^{12}$& 1.5$*10^{8}$\\
$10^{-16}$& $10^{-8}$&200023& 1.39$*10^{-14}$& 0.866& 3.90$*10^{-11}$& 5.3$*10^{11}$& 2.2$*10^{8}$\\
\end{tabular}
\vsn
The 5th column  verifies the statement of Theorem 1; all its entries should be at most 1. The 7th column shows how much the local error for equidistant 
mesh points exceeds that for the adaptive points used by ADMESH.  We see that a significant advantage of using adaption 
is observed for all values of $\e$ and $\delta$.  The gain grows when $\e$ or $\delta$ go to 0.
\vsn
{\bf Acknowledgments} $\;$ I thank L. Plaskota and P.  Przyby\l owicz for their comments on the manuscript.

\end{document}